\documentclass[a4paper,12pt]{article}
\usepackage{mathrsfs}
\usepackage{}
\usepackage[top=2.0cm,bottom=2.0cm,left=2.0cm,right=2.0cm]{geometry}
\usepackage{amssymb}
\usepackage{graphicx}
\usepackage{amsmath,amsthm,amssymb,lineno}
\usepackage{latexsym}
\usepackage{epstopdf}
\usepackage{setspace}
\usepackage{amsmath,color}
\usepackage{graphicx,booktabs,multirow}
\usepackage{latexsym, tabularx,shapepar}
\usepackage[all,2cell,dvips]{xy} \UseAllTwocells \SilentMatrices
\usepackage{appendix}
\usepackage{longtable}
\usepackage{cite}
\usepackage{CJK}
\usepackage{indentfirst}
\usepackage{array}
\usepackage[colorlinks,
           linkcolor=blue,
           anchorcolor=blue,
           citecolor=blue
           ]{hyperref}
\usepackage{amsthm}
\allowdisplaybreaks[4]
\graphicspath{{figures/}}

\newcommand \dbar {{\; \bar{} \hspace{-0.3em} \mathrm d}}

\newtheorem{theorem}{Theorem}[section]
\newtheorem{lemma}[theorem]{Lemma}
\newtheorem{remark}[theorem]{Remark}
\newtheorem{corollary}[theorem]{\rm\bfseries Corollary}

\newtheorem{conjecture}[theorem]{Conjecture}
\newtheorem{example}[theorem]{Example}
\newtheorem{problem}[theorem]{Problem}

\begin{document}
	
	\title{ On the eigenvalues and Seidel eigenvalues of chain graphs
	}
	\author{Zhuang Xiong\thanks{Corresponding author: zhuangxiong@hunnu.edu.cn}, Yaoping Hou \\
		{\small Key Laboratory of High Performance Computing and Stochastic Mathematics (Ministry of Education),} \\
		{\small College of Mathematics and Statistics, Hunan Normal University, Changsha, Hunan 410081,  China} \\
	}
	
	\date{}
	\maketitle
	\begin{abstract}
		
		In this paper we consider the eigenvalues and the Seidel eigenvalues of a chain graph. An$\dbar$eli\'{c}, da Fonseca, Simi\'{c}, and Du \cite{andelic2020tridiagonal} conjectured that there do not exist non-isomorphic cospectral chain graphs with respect to the adjacency spectrum. Here we disprove this conjecture. Furthermore, by considering the relation between the Seidel matrix and the adjacency matrix of a graph, we solve two problems on the number of distinct Seidel eigenvalues of a chain graph, which was posed by Mandal, Mehatari, and Das \cite{mandal2022spectrum}.
	\end{abstract}
	\baselineskip=0.25in
	
	\section{ Introduction}
\begin{figure}[htbp]\centering\hspace{0cm}
		\scalebox{0.6}{\includegraphics[width=8cm]{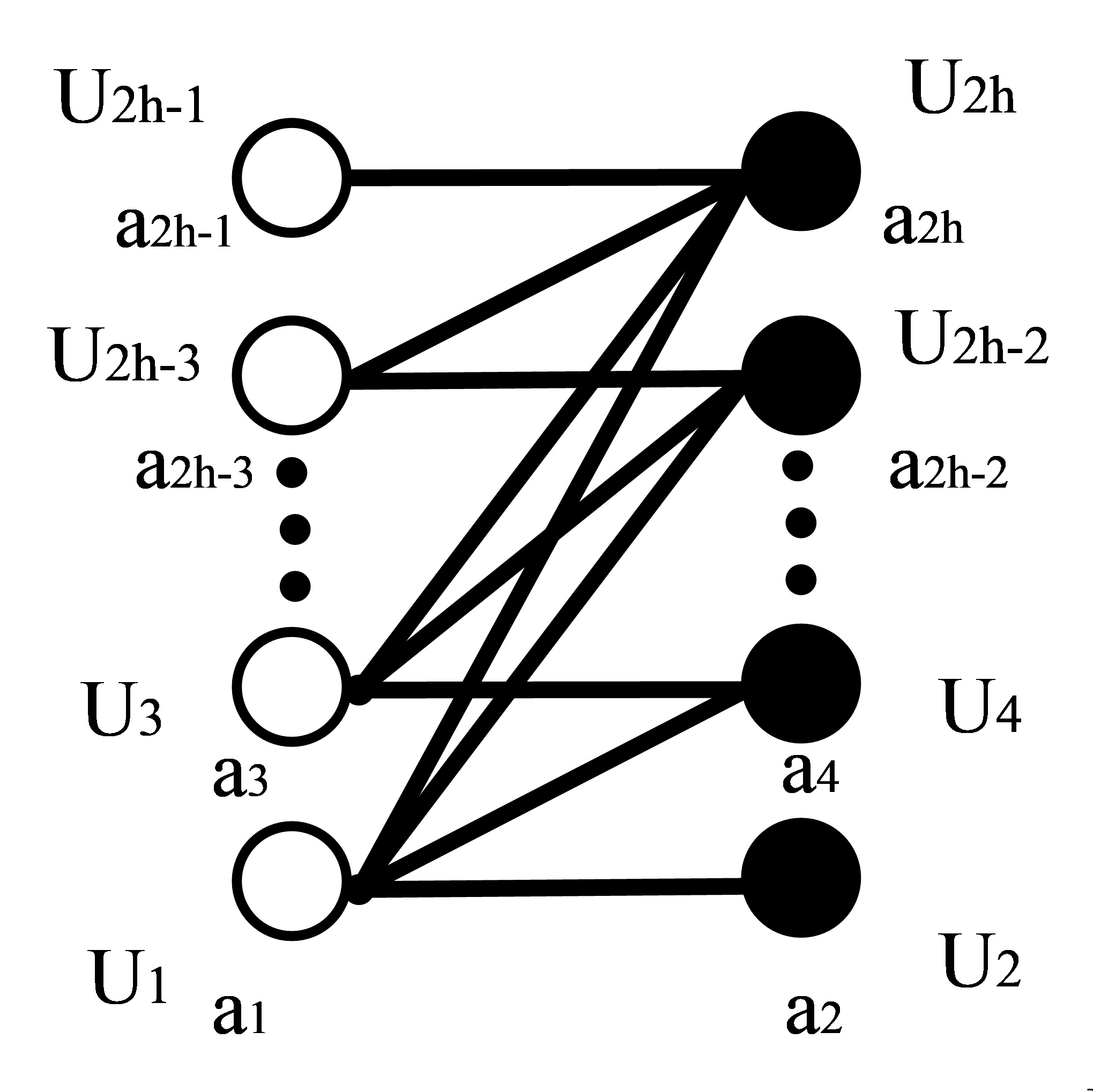}}\\
		~~~~~\caption{ The structure of a chain graph with binary string $0^{a_{1}}1^{a_{2}} \cdots 0^{a_{2h-1}}1^{a_{2h}}$, the thick line between distinct cells means that all vertices in one are adjacent to every vertex in the other. } ~~~~~~~~~~~~~~~~~~~~~~~~~~~~~~~~~~~~~~~~~~~~~~~~~~~~~~~~~~~~~~~~~~~~~~~~~~~~~~~~~~~~~~~~~~~~~~~~~~~~~~~~~~~~~~
		\label{fig_stru_chain_graph}
	\end{figure}

	Throughout this paper all graphs we consider are simple, finite, and undirected. Given a graph $G = (V(G), E(G))$ with vertex set $V(G)$ and edge set $E(G)$, its order and size are $|V(G)|$ and $ |E(G)|$, respectively. If two distinct vertices $v_i, v_j \in V(G)$ are adjacent in $G$, we write $v_i \sim v_j$, otherwise  $ v_i \nsim v_j$.  We denote by $A(G)$ and $S(G) = J - I - 2A(G)$ the adjacency matrix and the Seidel matrix of $G$, respectively, where $J$ is the all-ones matrix and $I$ is the identity matrix of order $n$. The eigenvalues of $A(G)$ and $S(G)$ are real, and we represent them as $\lambda_{1}(A) \geq \lambda_{2}(A) \geq \cdots \geq \lambda_{n}(A)$ and $\lambda_{1}(S) \geq \lambda_{2}(S) \geq \cdots \geq \lambda_{n}(S)$, respectively. Note that without causing confusion we abbreviate $A(G)$ and $S(G)$ as $A$ and $S$. In particular, when we mention the eigenvalues (resp., Seidel eigenvalues) and the spectrum (resp., Seidel spectrum) of $G$, we actually mean the eigenvalues and the spectrum of $A(G)$ (resp., $S(G)$). For more notions and notations about spectral graph theory, the reader is referred to \cite{brouwer2011spectra,cvetkovic2010introduction}.\\
	\indent A graph is called a chain graph if it is a $\{2K_2, C_3, C_5 \}$-free graph. It follows that chain graphs are bipartite. An equivalent definition is the following: a graph on $n$ vertices is a chain graph if it can be built by a binary string of length $n$. In other words, we can construct a chain graph by a recursive process. First there exists an isolated vertex (say white), then we can add an isolated white vertex or a vertex (say black) which is only adjacent to all existing white vertices at each step. We can see Fig. \ref{fig_stru_chain_graph} to understand this structure. Consequently, we represent a chain graph on $n$ vertices with a binary string $(a_1, a_2, \cdots, a_n)$ in which $a_i = 0$ if we add an isolated white vertex at step $i$, otherwise $a_i = 1$. This binary string is also written as $0^{a_1}1^{a_2}0^{a_3}1^{a_4} \cdots 0^{a_{2h-1}}1^{a_{2h}}$, where $0^{a_{k}}$ is composed of $a_{k}$ consecutive zeros and $1^{a_{l}}$ is composed of $a_{l}$ consecutive ones, for each $k = 1, \cdots, 2h-1$ and $l = 2, \cdots, 2h$.\\
	\indent Chain graphs play an important role in the class of connected bipartite graph in terms of their adjacency eigenvalues. In $2008$, Bell, Cvetkovi\'{c}, Rowlinson, and Simi\'{c} \cite{bell2008GraphsII} proved that the graphs which have the largest adjacency eigenvalue among the class of connected bipartite graphs are chain graphs, and named therein as double nested graphs. This result also discovered independently by Bhattacharya, Shmuel, and Peled in \cite{bhattacharya2008on}. After that, there has been a lot of researches on the mathematical properties of chain graphs. In particular, the researches of related matrices have received wide attention in recent years. In \cite{andelic2011on}, An$\dbar$eli\'{c}, da Fonseca, Simi\'{c}, and To\v{s}i\'{c} studied the adjacency matrix and gave some lower and upper bounds for the largest adjacency eigenvalue of double nested graphs. In $2021$, the distance matrix of a chain graph are considered by Alazemi,  An$\dbar$eli\'{c}, Koledin, and Stani\'{c} \cite{alazemi2021eigenvalue}.
	They obtained eigenvalue-free interval of distance matrices of threshold and chain graphs. In $2022$, Alazemi, An$\dbar$eli\'{c}, Das, and da Fonseca \cite{alazemi2022chain} gave some properties of the spectrum of the Laplacian matrix of a chain graph in terms of its degree sequence. Very recently, the spectral properties of the Seidel matrix of a chain graph was investigated by Mandal, Mehatari, and Das \cite{mandal2022spectrum}.\\
	\indent In this work, we deal with the eigenvalues of the adjacency matrices and the Seidel matrices of chain graphs. In Section \ref{sec:pre}, we introduce some basic notations and results which will be used in the sequel. In Section \ref{sec:eig}, we focus attention on the eigenvalues of chain graphs. We determine the chain graphs with distinct eigenvalues and disprove the conjecture posed in \cite{andelic2020tridiagonal} which say that there do not exist a pair of non-isomorphic cospectral connected chain graphs. The Seidel eigenvalues of chain graphs are considered in Section \ref{sec:Seideleig}. By the same motivation as the work of Xiong and Hou in \cite{zhuang2022eigenvalue},  we characterize the distribution of Seidel eigenvalues of chain graphs, which solve two problems proposed in \cite{mandal2022spectrum} on the numbers of positive and negative Seidel eigenvalues. Finally, by using the results obtained previously we determine the chain graphs with distinct Seidel eigenvalues.\\
	\indent

	\section{Preliminaries}\label{sec:pre}
	
	We fix notation and introduce several important results which will be used in the sequel.
	Given a square matrix $Q$, the inertia of $Q$ is a triple $(n_{+}(Q), n_{0}(Q), n_{-}(Q))$, where $n_{0}(Q)$ is the nullity of $Q$ and $n_{+}(Q)$ (resp., $n_{-}(Q)$) is the number of positive (resp., negative) eigenvalues of $Q$. The spectrum of $Q$ is denoted by $Spec(Q) = \{[\lambda_1]^{m_1}, \cdots, [\lambda_k]^{m_k}\}$, where $\lambda_i$, for each $i =1, \cdots ,k$, is an eigenvalue of $Q$ with multiplicity $m_{i}$ and the exponents indicate multiplicities. The characteristic polynomial of $Q$ is written as $P(Q, x) = \det(xI - Q)$. We denote by rank$(Q)$ the rank of $Q$.\\
	\indent Let $G$ be a chain graph built by the binary string $0^{a_{1}}1^{a_{2}} \cdots 0^{a_{2h-1}}1^{a_{2h}}$. We denote by \textbf{$\Pi$} = $ U_{1} \cup U_{2} \cup \cdots \cup U_{2h-1} \cup U_{2h}$ an equitable partition of the vertex set $V(G)$. The vertex subset $U_{k}$ forms a coclique, where $|U_{k}| = a_{k}$, $1 \leq k \leq 2h $. Therefore, the quotient matrices $\widetilde{A}(G)$ and $\widetilde{S}(G)$ of the adjacency matrix $A(G)$ and the Seidel matrix $S(G)$ with respect to partition \textbf{$\Pi$} are
	\begin{equation*}
		\begin{array}{cc}
			\begin{pmatrix}
				0 & a_2 & 0 & a_4 & \cdots & 0 & a_{2h}\\
				a_1 & 0 & 0 & 0 & \cdots & 0 & 0\\
				0 & 0 & 0 & a_4 & \cdots & 0 & a_{2h}\\
				a_1 & 0 & a_3 & 0 & \cdots & 0 & 0\\
				\vdots & \vdots & \vdots & \vdots & \ddots & \vdots & \vdots\\
				0 & 0 & 0 & 0 & \cdots & 0 & a_{2h}\\
				a_1 & 0 & a_3 & 0 & \cdots & a_{2h-1} & 0
			\end{pmatrix} and
			\begin{pmatrix}
				a_1-1 & -a_2 & a_3 & -a_4 & \cdots & a_{2h-1} & -a_{2h}\\
				-a_1 & a_2-1 & a_3 & a_4 & \cdots & a_{2h-1} & a_{2h}\\
				a_1 & a_2 & a_3-1 & -a_4 & \cdots & a_{2h-1} & -a_{2h}\\
				-a_1 & a_2 & -a_3 & a_4-1 & \cdots & a_{2h-1} & a_{2h}\\
				\vdots & \vdots & \vdots & \vdots & \ddots & \vdots & \vdots\\
				a_1 & a_2 & a_3 & a_4 & \cdots & a_{2h-1}-1 & -a_{2h}\\
				-a_1 & a_2 & -a_3 & a_4 & \cdots & -a_{2h-1} & a_{2h}-1
			\end{pmatrix},
		\end{array}
	\end{equation*}respectively. In the rest of this paper, $\widetilde{A}(G)$ and $\widetilde{S}(G)$ always denote the two quotient matrix given here. For a vertex $v \in V(G)$, let $N_{G}(v)$ denote the vertex set in which the vertices are adjacent to $v$. When the text is clear, we will omit the subscript $G$ and two vertices $u$ and $v$ are called duplicates if $N(u) = N(v)$.  It follows that, for the partition \textbf{$\Pi$}, each pair of vertices which are in the same cell are a pair of duplicate vertices.
	
	\indent An important result in spectral graph theory is the Interlacing Theorem. We record it as follows.
	\begin{lemma}\cite[Corollary 1.3.12]{cvetkovic2010introduction}\label{lem:interlacing}
		Let $G$ be a graph with $n$ vertices and eigenvalues $\lambda_1 \geq \lambda_2 \geq \cdots \geq \lambda_n$, and let $H$ be an induced subgraph of $G$ with $m$ vertices. If the eigenvalues of $H$ are $\mu_1 \geq \mu_2 \geq \cdots \geq \mu_m$ then $\lambda_{n-m+i} \leq \mu_i \leq \lambda_i$ $(i = 1, \cdots, m).$
	\end{lemma}
	We end up this section by introducing the Courant-Weyl inequalities.
	\begin{lemma}\cite[Theorem 1.3.15]{cvetkovic2010introduction}\label{lem:cw}
		Let $M$ and $N$ be $n \times n$ Hermitian matrices. Then
		\begin{gather*}
			\lambda_{i}(M + N) \leq \lambda_{j}(M) + \lambda_{i-j+1}(N) \qquad (1 \leq j \leq i \leq n), \\
			\lambda_{i}(M + N) \geq \lambda_{j}(M) + \lambda_{i-j+n}(N) \qquad (1 \leq i \leq j \leq n),
		\end{gather*}
		where the eigenvalues of $M + N$, $M$, and $N$ are all arranged in non-increasing order.
	\end{lemma}
	
	\section{ The eigenvalues of chain graphs}\label{sec:eig}
	In this section, we mainly consider spectral properties of the adjacency matrix $A(G)$ for a chain graph $G$. We start by introducing some results on the eigenvalues of a chain graph.
	
	\begin{lemma}\cite[Theorem 4.3]{andelic2015some}\label{lem:spectrumA}
		The spectrum of a connected chain graph consists of $2h$ distinct nonzero eigenvalues (determined by the quotient matrix $\widetilde{A}(G)$) and of $0$ with multiplicity $n-2h$.
	\end{lemma}
	
	\begin{lemma}\label{lem:efree}\cite[Theorem 3.3]{andelic2016location}
		A chain graph $G$ has no eigenvalue in the interval $[-\frac{1}{2}, 0) \cup (0, \frac{1}{2}]$.
	\end{lemma}
	
	The problem: ``Which graphs have distinct eigenvalues?'' was proposed by Harary and Schwenk \cite{harary1974which} in 1974. Next we will determine all chain graphs with distinct eigenvalues.
	
	\begin{theorem}
		Let $G$ be a chain graph with binary $0^{a_{1}}1^{a_{2}} \cdots 0^{a_{2h-1}}1^{a_{2h}}$. Then all eigenvalues of $G$ are distinct if and only if one of the followings holds. \\
		\indent (\romannumeral1) $a_{1} = \cdots = a_{2h} = 1$;\\
		\indent (\romannumeral2) $a_{i} = 2$ and $a_{j} = 1$ for $j \neq i$.
	\end{theorem}
	\begin{proof}
		In the case of $(\romannumeral1)$, we know that $\widetilde{A}(G) = A(G)$ and the result follows easily from Lemma \ref{lem:spectrumA}.
		In the case of $(\romannumeral2)$, using Lemma \ref{lem:spectrumA}, we obtain that the spectrum of $A(G)$ consists of $0$ with multiplicity $1$ and the $2h$ distinct eigenvalues of $\widetilde{A}(G)$. Thus the sufficiency follows and we show the necessity below.\\
		\indent  Since all the eigenvalues of $A(G)$ are distinct, in particular, the multiplicity of eigenvalue $0$ is at most $1$. So we have $n-2h \leq 1$ by considering the spectrum of $A(G)$. The case of $n-2h = 0$ leads to (\romannumeral1) and the case of $n-2h = 1$ leads to (\romannumeral2), respectively.
	\end{proof}
	
	In \cite{andelic2020tridiagonal} An$\dbar$eli\'{c}, da Fonseca, Simi\'{c}, and Du conjectured that there exist no a pair of non-isomorphic cospectral connected chain graphs and we state this conjecture as follows.  Note that the multiplicity of the eigenvalue  $0$ of a chain graph is $n-2h$. Hence a pair of cospectral chain graphs must have the same parameters $n$ and $h$. We will disprove the conjecture by considering the chain graphs with the parameter $h$ equal to $2$ in Theorem \ref{tho:h2¡ª¡ªyes}. In spite of that, a weakly version of the conjecture is true: In Theorem \ref{tho:h1_no} it will be shown that for a chain graph $G$ with binary string $0^{a_1} 1^{a_2}$, there does not exist a connected chain graph which is cospectral but non-isomorphic to $G$, i.e., the conjecture holds for the class of chain graphs  with $h$ equal to $1$.
	\begin{conjecture}\cite[Conjecture 3.3]{andelic2020tridiagonal}\label{conj:noco}
		In a class of connected chain graphs there do not exist nonisomorphic cospectral chain graphs with respect to the adjacency spectrum.
	\end{conjecture}
	
	In the class of chain graphs on $n$ vertices with $h = 1$, for any chain graph with binary $0^{a_1} 1^{a_2}$, say, $G$, the spectrum of $G$ is $\{\sqrt{a_1a_2}, [0]^{n-2}, -\sqrt{a_1a_2} \}$. If there exists another chain graph $H$ with binary $0^{b_1} 1^{b_2}$, which is cospectral with $G$, then we have $b_1 b_2 = a_1 a_2$ and $b_1 + b_2 = n = a_1 + a_2$. It follows that $a_1 = b_1$ and $a_2 = b_2$ or $a_1 = b_2$ and $a_2 = b_1$, which means that $G$ is isomorphic to $H$. We record this result as the following theorem.
	\begin{theorem}\label{tho:h1_no}
		Given a chain graph $G$ on $n$ vertices with binary string $0^{a_{1}}1^{a_{2}}$, there does not exist a chain graph which is non-isomorphic but cospectral to $G$.
	\end{theorem}

	In the next theorem, we construct infinitely many pairs of cospectral chain graphs for $h = 2$. This provide a counterexample to disprove Conjecture \ref{conj:noco}.
	
	\begin{theorem}\label{tho:h2¡ª¡ªyes}
		Two chain graphs $G$ and $H$ with binary strings $0^{a_{1}}1^{a_{2}}0^{a_{3}}1^{a_{4}}$ and $0^{a_{2}}1^{a_{1}}0^{a_{4}}1^{a_{3}}$, respectively, are cospectral when $a_{1}a_{4} = a_{2}a_{3}$. Furthermore, if $a_{1} \neq a_{2}$ and $a_{1} \neq a_{3}$, then $G$ is non-isomorphic to $H$.
	\end{theorem}
	\begin{proof}As we mentioned above, both $G$ and $H$ have the eigenvalue $0$ with multiplicity $n-4$. To investigate the remaining $4$ eigenvalues, we shall consider the eigenvalues of related quotient matrices. By a simple computation, we obtain that characteristic polynomials of $\widetilde{A}(G)$ and $\widetilde{A}(H)$ are $x^4 - (a_{1}a_{2} + a_{1}a_{4} + a_{3}a_{4})x^2  + a_{1}a_{2}a_{3}a_{4}$ and $x^4 - (a_{1}a_{2} + a_{2}a_{3} + a_{3}a_{4})x^2  + a_{1}a_{2}a_{3}a_{4}$, respectively. If $a_{1}a_{4} = a_{2}a_{3}$, then the two polynomials are equivalent, that is, $\widetilde{A}(G)$ and $\widetilde{A}(H)$ share the same spectrum. So far we have proved the first half of the theorem. In what follows, we will show that, under our restrictions, $G$ is non-isomorphic to $H$.\\
		\indent  Suppose to the contrary that $G$ is isomorphic to $H$ when $a_{1} \neq a_{2}$ and $a_{1} \neq a_{3}$. We observe that the degree sequences of $G$ and $H$ are $((a_{1})^{a_{2}}, ({a_{2}+a_{4}})^{a_{1}}, ({a_{1}+a_{3}})^{a_{4}}, (a_{4})^{a_{3}})$ and $((a_{2})^{a_{1}}, ({a_{2}+a_{4}})^{a_{3}}, ({a_{1}+a_{3}})^{a_{2}}, (a_{3})^{a_{4}})$, respectively, where the superscripts denote the numbers of vertices with corresponding degrees. Note that $a_{i} \geq 1$ for each $i = 1,2,3,4$. By isomorphism we know that $G$ and $H$ share the same degree sequence. Now, we first consider the numbers of vertices with degree ${a_{2}+a_{4}}$ of $G$ and $H$. We have seen that there have $a_{1}$ (resp., $a_{3}$) vertices with degree ${a_{2}+a_{4}}$ in $G$ (resp., $H$) so far. Since $a_{1} \neq a_{3}$, so there exists other vertices with degree ${a_{2}+a_{4}}$ and we shall distinguish three cases in terms of this. \\ 
		\indent Case $1$: $a_{1} + a_{3} = a_{2} + a_{4}$. Then $a_{1} + a_{4} = a_{3} + a_{2}$ by considering the numbers of vertices with degree $a_{1} + a_{3}$ in $G$ and $H$. Because $a_{1} \neq a_{3}$, we have $a_{2} \neq a_{4}$. Next we consider the numbers of the vertices with degree $a_{1}$ in the two graphs. Since $a_{1} \neq a_{1} + a_{3} = a_{2} + a_{4}$, $a_{1} \neq a_{2}$, and $a_{1} \neq a_{3}$, there does not exist a vertex with degree $a_1$ in $H$, a contradiction.\\
		\indent Case 2: $a_{1} = a_{2} + a_{4}$. Since $a_{1} \neq a_{1} + a_{3}$, $a_{2} + a_{4} \neq a_{4}$, and $a_{1} \neq a_{3}$, we have $a_{1} + a_{2} = a_{3}$ by comparing the numbers of the vertices with degree $a_{2} + a_{4}$ in $G$ and $H$. By considering the numbers of the vertices with degree $a_{1} + a_{3}$ in two graphs we have $a_{2} = a_{4}$. Finally, we compare the numbers of the vertices with degree $a_{2}$ in the two graphs. Since $a_{1} \neq a_{3}$, there exist other vertices with degree $a_{2}$ in $H$, it follows that $a_{2} = a_{3}$. Recall that $a_{1} + a_{2} = a_{3}$, which leads to $a_{1} = 0$, a contradiction.\\
		\indent Case 3: $a_{3} = a_{2} + a_{4}$. Since $a_1 \neq a_2$, $a_1 \neq a_1 + a_3$, and $a_1 \neq a_3$, then there exist $a_2$ vertices with degree $a_1$ in $G$, but there do not exist such vertices in $H$, which leads to a contradiction. 
	\end{proof}
	
	We illustrate the theorem above with an example as the ending of this section.
	\begin{example}\label{ex:h2} The chain graph $G$ with binary string $0^1 1^2 0^2 1^4$ and the chain graph $H$ with binary string $0^2 1^1 0^4 1^2$, are a pair of cospectral non-isomorphic chain graphs. The two chain graphs above share the spectrum $\{ -3.57, -1.12, [0]^5, 1.12, 3.57 \}$ but the degree sequences of $G$ and $H$ are $(1, 1, 3, 3, 3, 3, 4, 4, 6)$ and $(2, 2, 2, 2, 2, 3, 3, 6, 6)$, respectively.
	\end{example}
	
	\begin{figure}[htbp]\centering\hspace{0cm}
		\scalebox{0.6}{\includegraphics[width=10cm]{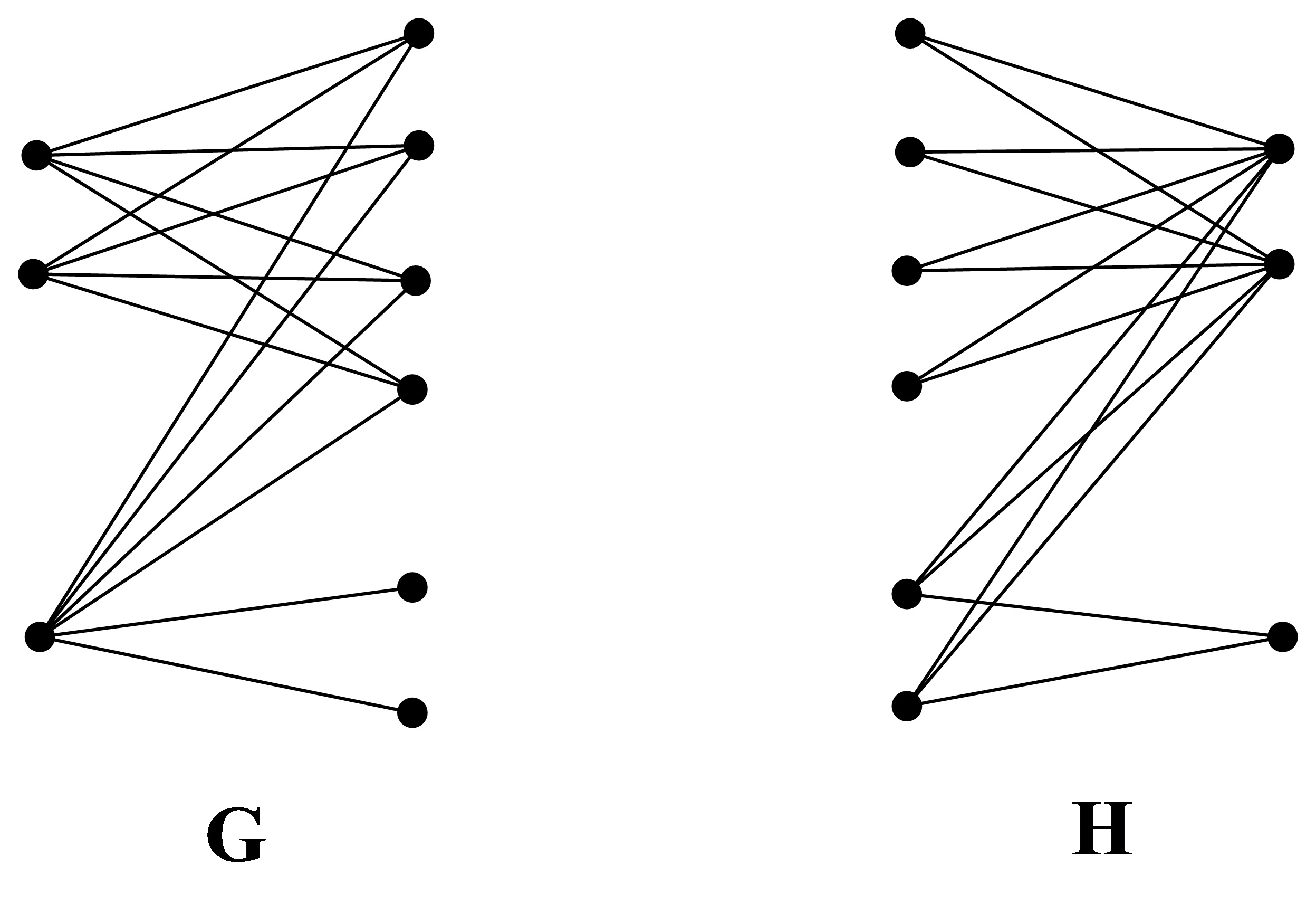}}\\
		\caption{ The two chain graphs mentioned in Example \ref{ex:h2}}
		\label{fig:h2_ex}
	\end{figure}

	\section{ The Seidel eigenvalues of chain graphs}\label{sec:Seideleig}
	We will now turn our attention to the Seidel eigenvalues of a chain graph. Similar to the method in Section \ref{sec:eig}, we begin with considering the quotient matrix of a Seidel matrix for a chain graph, from which we obtain some characterizations of the Seidel eigenvalues of a chain graph. Our results also solve some problems proposed in \cite{mandal2022spectrum} on the number of distinct Seidel eigenvalues and the inertia of the Seidel matrix for a chain graph. We start by introducing a lemma about the Seidel eigenvalue $-1$ of a connected graph according to the size of its independent set.
	
	\begin{lemma}\label{lem:seidel_-1} Let $G = (V, E)$ be a connected graph on $n$ vertices and $W$ a subset of $V(G)$ with size $l$, $2 \leq l \leq n$. If $W$ is an independent set and $N(u) = N(v)$ for all $u, v \in V(W)$, then $-1$ is an Seidel eigenvalues of $G$ with multiplicity at least $l-1$.
	\end{lemma}
	\begin{proof} Let $W = \{ w_{1}, \cdots, w_{l} \} \subseteq V(G)$ is an independent set of $G$ such that $N(w_{i})= N(w_{j}) = C = \{v_{1}, \cdots, v_{k}\}$ holds for $1 \leq i, j \leq l$. Set $D = V(G) \setminus (W \cup C) = \{ u_{1}, \cdots u_{n-k-l} \}$. Then $V(G) = W \cup C \cup D$ is a vertex partition of $G$ for which the induced graph of $G$ on $W$ is an empty graph of order $l$. Moreover, there exist no edges between $W$ and $D$ in $G$, so the submatrix $S(W,D)$ of the Seidel matrix of $G$ induced on the row set $W$ and the column set $D$  equals to $J_{l \times (n-l-k)}$. Similarly we have $S(W,C) = -J_{l \times k}$, where $S(W,C)$ is the submatrix of $S(G)$ induced on the row set $W$ and the column set $C$. Therefore, the Seidel matrix $S(G)$ can be denoted as
		$$
		\bordermatrix{
			& W & D & C \cr
			W & (J-I)_{l \times l} & J_{l \times (n-l-k)} & -J_{l \times k} \cr
			D & J_{(n-k-l) \times l} & E & F \cr
			C & -J_{k \times l} & F^{T} & H \cr
		},
		$$
		where the block matrix $E = S(D,D)$ is the submatrix of $S(G)$ induced on the vertex set $D$ and $F^T$ denote the transpose of the matrix $F$. By the same way we have $F = S(D,C)$ and $H = S(C,C)$. Take $x^{(p)} \in \mathbb{R}^{n} (2 \leq p \leq l) $ is a vector defined on the $V(G)$ satisfying that $x ^{(p)}(v_{1}) = 1, x ^{(p)}(v_{p}) = -1$ and $x ^{(p)}(v) = 0$ for $v \notin \{ v_{1}, v_{p} \}$. It follows that $S(G) \cdot x^{(p)} = -x^{(p)}$ and $x^{(2)}, \cdots, x^{(l)}$ are linearly independent. Therefore, $-1$ is a Seidel eigenvalue of $G$ with multiplicity at least $l-1$.
	\end{proof}
	
	Because each cell of partition \textbf{$\Pi$} of a given chain graph is an independent set satisfying the conditions of Lemma \ref{lem:seidel_-1}, we can directly obtain the following corollary.
	\begin{corollary} Let $G$ be a chain graph with binary string $0^{a_{1}}1^{a_{2}} \cdots 0^{a_{2h-1}}1^{a_{2h}}$. Then $-1$ is a Seidel eigenvalue with multiplicity at least $n-2h$.
	\end{corollary}
	
	Follow the above notations and note that the Seidel eigenvalues of $G$ can be divided into two types according to their associated eigenvectors. First we know that all $2h$ eigenvalues of $\widetilde{S}(G)$ are the eigenvalues of $S(G)$.  These eigenvalues are of the type $1$: each of the eigenvectors of $S(G)$ corresponding to an eigenvalue of the quotient matrix $\widetilde{S}(G)$ is constant on the parts of the partition \textbf{$\Pi$}, that is, the entries corresponding to the vertices which belong to the same cell are equal. Using the result of Lemma \ref{lem:seidel_-1}, we define the $n-2h$ linearly independent eigenvectors associated to the eigenvalue $-1$ of type $2$, which are all orthogonal to the eigenvectors corresponding to the eigenvalues of type $1$. From these observations we have that the characteristic polynomial of $S(G)$ is $P(S(G),x) = (x+1)^{n-2h} \cdot P(\widetilde{S}(G),x) = (x+1)^{n-2h} \cdot \det(\widetilde{S}(G) - xI_{2h})$. To obtain more information on the eigenvalues of $S(G)$ we will apply a series of elementary row operations to the matrix $\widetilde{S}(G) - xI_{2h}$. First, performing the following elementary row transformations:\\
	\indent $ R_{2k-1} \leftarrow R_{2k-1} - R_{2k+1},\quad k = 1, 2, \cdots, h-1$,\\
	\indent $ R_{2k} \leftarrow R_{2k} - R_{2k-2},\quad k = 2, \cdots, h$,\\
	we have
	\[ \begin{pmatrix}
		-1-x & -2a_2 & 1+x & 0 & \cdots & 0 & 0\\
		-a_1 & a_2-1-x & a_3 & a_4 & \cdots & a_{2h-1} & a_{2h}\\
		0 & 0 & -1-x & -2a_4 & \cdots & 0 & 0\\
		0 & 1+x & -2a_3 & -1-x & \cdots & 0 & 0\\
		\vdots & \vdots & \vdots & \vdots & \ddots & \vdots & \vdots\\
		a_1 & a_2 & a_3 & a_4 & \cdots & a_{2h-1}-1-x & -a_{2h}\\
		0 & 0 & 0 & 0 & \cdots & -2a_{2h-1} & -1-x
	\end{pmatrix}. \]
	Continuing with $R_{2k-1} \leftrightarrow R_{2k}$, i.e., interchanging the positions of row $2k-1$ and row $2k$, for each $k = 1,2, \cdots, h,$ we derive\\
	\[ \begin{pmatrix}
		-a_1 & a_2-1-x & a_3 & a_4 & \cdots & a_{2h-1} & a_{2h}\\
		-1-x & -2a_2 & 1+x & 0 & \cdots & 0 & 0\\
		0 & 1+x & -2a_3 & -1-x & \cdots & 0 & 0\\
		0 & 0 & -1-x & -2a_4 & \cdots & 0 & 0\\
		\vdots & \vdots & \vdots & \vdots & \ddots & \vdots & \vdots\\
		0 & 0 & 0 & 0 & \cdots & -2a_{2h-1} & -1-x\\
		a_1 & a_2 & a_3 & a_4 & \cdots & a_{2h-1}-1-x & -a_{2h}
	\end{pmatrix}. \]
	Then by\\
	\indent $  R_{1} \leftarrow \sum_{k = 1}^h R_{2k-1}$,\\
	\indent $  R_{2h} \leftarrow \sum_{k = 1}^h R_{2k}$,\\
	we obtain
	\[ \begin{pmatrix}
		-a_1 & a_2 & -a_3 & a_4 & \cdots & -a_{2h-1} & a_{2h}-1-x\\
		-1-x & -2a_2 & 1+x & 0 & \cdots & 0 & 0\\
		0 & 1+x & -2a_3 & -1-x & \cdots & 0 & 0\\
		0 & 0 & -1-x & -2a_4 & \cdots & 0 & 0\\
		\vdots & \vdots & \vdots & \vdots & \ddots & \vdots & \vdots\\
		0 & 0 & 0 & 0 & \cdots & -2a_{2h-1} & -1-x\\
		a_1-1-x & -a_2 & a_3 & -a_4 & \cdots & a_{2h-1} & -a_{2h}
	\end{pmatrix}. \]
	By\\
	\indent $  R_{1} \leftarrow R_{1} - \sum_{k = 2}^h \frac{1}{2} R_{2k-1}$,\\
	\indent $  R_{2h} \leftarrow R_{2h} + \sum_{k = 2}^h \frac{1}{2} R_{2k-1}$,\\
	we have
	\[ \begin{pmatrix}
		-a_1 & a_2-\frac{1}{2}(1+x) & 0 & a_4 & \cdots & 0 & a_{2h}-\frac{1}{2}(1+x)\\
		-1-x & -2a_2 & 1+x & 0 & \cdots & 0 & 0\\
		0 & 1+x & -2a_3 & -1-x & \cdots & 0 & 0\\
		0 & 0 & -1-x & -2a_4 & \cdots & 0 & 0\\
		\vdots & \vdots & \vdots & \vdots & \ddots & \vdots & \vdots\\
		0 & 0 & 0 & 0 & \cdots & -2a_{2h-1} & -1-x\\
		a_1-1-x & -a_2+\frac{1}{2}(1+x) & 0 & -a_4 & \cdots & 0 & -a_{2h}-\frac{1}{2}(1+x)
	\end{pmatrix}. \]
	Proceeding with\\
	\indent $  R_{1} \leftarrow R_{1} + R_{2h}$,\\
	\indent $  R_{2h} \leftarrow R_{2h} - \sum_{k = 1}^{h-1} \frac{1}{2} R_{2k}$,\\
	we derive\\
	\[ \begin{pmatrix}
		-1-x & 0 & 0 &0 & \cdots & 0 & -1-x\\
		-1-x & -2a_2 & 1+x & 0 & \cdots & 0 & 0\\
		0 & 1+x & -2a_3 & -1-x & \cdots & 0 & 0\\
		0 & 0 & -1-x & -2a_4 & \cdots & 0 & 0\\
		\vdots & \vdots & \vdots & \vdots & \ddots & \vdots & \vdots\\
		0 & 0 & 0 & 0 & \cdots & -2a_{2h-1} & -1-x\\
		a_1-\frac{1}{2}(1+x) & \frac{1}{2}(1+x) & 0 & 0 & \cdots & -\frac{1}{2}(1+x) & -a_{2h}-\frac{1}{2}(1+x)
	\end{pmatrix}. \]
	Finally, using\\
	\indent $  R_{2h} \leftarrow 2R_{2h}$,\\
	\indent $  R_{2h} \leftarrow R_{2h} - R_{1}$,\\
	\indent $  R_{k} \leftarrow -R_{k}, k = 1, \cdots, 2h$,\\
	we obtain a matrix
	\[ F_{x} = \begin{pmatrix}\tag{4.1}\label{equ:matrix_fin}
		1+x & 0 & 0 &0 & \cdots & 0 & 1+x\\
		1+x & 2a_2 & -1-x & 0 & \cdots & 0 & 0\\
		0 & -1-x & 2a_3 & 1+x & \cdots & 0 & 0\\
		0 & 0 & 1+x & 2a_4 & \cdots & 0 & 0\\
		\vdots & \vdots & \vdots & \vdots & \ddots & \vdots & \vdots\\
		0 & 0 & 0 & 0 & \cdots & 2a_{2h-1} & 1+x\\
		-2a_1 & -1-x & 0 & 0 & \cdots & 1+x & 2a_{2h}
	\end{pmatrix}. \]
	\indent By considering all the row operations above, we have $\det(\widetilde{S}(G) - xI_{2h}) = (-1)^{h} \cdot \frac{1}{2} \det(F_{x})$, which means that the roots of the equation $\det(F_{x}) = 0$ are the eigenvalues of $\widetilde{S}(G)$. Thus we can derive the following lemmas by some simple row and column operations on some matrices we need. Note that Lemmas \ref{lem:q_e_-1} and \ref{lem:q_e_oth} are also obtained by \cite[Theorem 2.2]{mandal2022spectrum} and \cite[Theorem 2.3]{mandal2022spectrum}, respectively. Here we also give our proofs for self contained.
	
	\begin{lemma}\label{lem:q_e_-1}
		$-1$ is a simple eigenvalue of the quotient matrix $\widetilde{S}(G)$.
	\end{lemma}
	\begin{proof}
		Replacing $x$ with $-1$ from (\ref{equ:matrix_fin}), we have the matrix
		\[ F_{-1} = \begin{pmatrix}
			0 & 0 & 0 &0 & \cdots & 0 & 0\\
			0 & 2a_2 & 0 & 0 & \cdots & 0 & 0\\
			0 & 0 & 2a_3 & 0 & \cdots & 0 & 0\\
			0 & 0 & 0 & 2a_4 & \cdots & 0 & 0\\
			\vdots & \vdots & \vdots & \vdots & \ddots & \vdots & \vdots\\
			0 & 0 & 0 & 0 & \cdots & 2a_{2h-1} & 0\\
			-2a_1 & 0 & 0 & 0 & \cdots & 0 & 2a_{2h}
		\end{pmatrix}. \]
		Since $\det(F_{-1}) = 0$ we know that $-1$ is an eigenvalue of $\widetilde{S}(G)$. Moreover, because the submatrix of $F_{-1}$ obtained by deleting the first row and the first column is a tridiagonal matrix with positive diagonal entries, we have that rank$(\widetilde{S}(G) + I_{2h}) = $ rank$ (F_{-1}) = 2h - 1$.  So $-1$ is a simple eigenvalue of $\widetilde{S}(G)$.
	\end{proof}
	
	We next consider the multiplicities of the remaining eigenvalues of $\widetilde{S}(G)$.
	\begin{lemma}\label{lem:q_e_oth}
		If $\lambda \neq 1$ is an eigenvalue of $\widetilde{S}(G)$, then the multiplicity of $\lambda$ is at most $2$.
	\end{lemma}
	\begin{proof}
		Since $\lambda$ is an eigenvalue of $\widetilde{S}(G)$, we have $\det(F_{\lambda}) = 0$. Deleting the first and the last row and the last two columns, we obtain a matrix
		\[ F_{\lambda}^{(1)} = \begin{pmatrix}
			1+\lambda & 2a_{2} & -1-\lambda &0 & \cdots & 0 & 0\\
			0 & -1-\lambda & 2a_{3} & 1+\lambda & \cdots & 0 & 0\\
			0 & 0 & 1+\lambda & 2a_{4} & \cdots & 0 & 0\\
			\vdots & \vdots & \vdots & \vdots & \ddots & \vdots & \vdots\\
			0 & 0 & 0 & 0 & \cdots & 1+\lambda & 2a_{2h-2}\\
			0 & 0 & 0 & 0 & \cdots & 0 & -1-\lambda
		\end{pmatrix}. \]
		Because $F_{\lambda}^{(1)}$ is an upper triangular matrix with nonzero diagonal entries, we have rank$(\widetilde{S}(G) - \lambda I_{2h}) = $ rank$ (F_{\lambda}) \geq 2h - 2$. So the multiplicity of $\lambda$ is at most $2$.
	\end{proof}
	
	\begin{remark}\label{rem:m_s} We have known from the previous section that, for any chain graph $G$, the eigenvalues of $\widetilde{A}(G)$ are simple. This does not hold for $\widetilde{S}(G)$. For example, the spectrum of $\widetilde{S}(G)$ of a chain graph $G$ with binary string $0^1 1^2 0^2 1^2 0^2 1^1$ is $\{ [-3.4721]^{2}, -1, 1, [5.4721]^{2} \}$. Given a chain graph with binary string $0^{a_{1}}1^{a_{2}} \cdots 0^{a_{2h-1}}1^{a_{2h}}$, if we denote by $M_{S}(G)$ the number of distinct Seidel eigenvalues of $G$, it was shown in \cite[Corollary 2.7]{mandal2022spectrum} that $h+1 \leq M_{S}(G) \leq 2h$. In \cite{mandal2022spectrum}, a problem on this parameter was proposed by Mandal, Mehatari, and Das, which will be described as follows and we will give an answer of it.
	\end{remark}
	
	\begin{problem}\cite[Problem 2.9]{mandal2022spectrum}\label{pro:ext_m_s} Is there exist a chain graph for which $ h+1 < M_{S}(G) < 2h$?
	\end{problem}
	
	Note that for $h = 2$ there exists no chain graph satisfying the condition of the above problem. Thus we consider the case of $h = 3$ and we want to find a class of chain graphs which meet this condition. Unfortunately, we found only a few sporadic examples, and the chain graphs of these examples do not seem to fit into a class based on their parameters. We list below two chain graphs meeting the requirements, which positively address Problem \ref{pro:ext_m_s}.
	\begin{example}\label{ex:solve_m_s}
		For $h = 3$, the spectra of chain graph $G_{1}$ with binary string $0^2 1^4 0^2 1^6 0^2 1^2$ and chain graph $G_{2}$ with binary string $0^5 1^2 0^4 1^4 0^2 1^1$
		are $\{21.1168, 11, 3.8831, [-1]^{13}, [-7]^2 \}$ and $\{13.7445, 7, 2.2554, [-1]^{13}, [-5]^2 \}$, respectively.
	\end{example}
	
	In Section \ref{sec:eig}, we determined the chain graphs with distinct eigenvalues. By the same motivation we will characterize the chain graphs with distinct Seidel eigenvalues. First we introduce the following result on the determinant of a special tridiagongal matrix $M_{k}(c)$, where
	\[ M_k(c) = \begin{pmatrix}
		c & 1 & 0 & 0 & \cdots & 0 & 0\\
		1 & c & 1 & 0 & \cdots & 0 & 0\\
		0 & 1 & c & 1 & \cdots & 0 & 0\\
		\vdots & \vdots & \vdots & \vdots & \ddots & \vdots & \vdots\\
		0 & 0 & 0 & 0 & \cdots & c & 1\\
		0 & 0 & 0 & 0 & \cdots & 1 & c
	\end{pmatrix}, ~~~~~~ c \in \mathbb{C}, ~~~k \in \mathbb{N}. \]
	
	In the following lemma we denote the determinant $\det(M_k(c))$ of the $k \times k$ tridiagonal matrix $M_k(c)$ by $D_k(c)$.
	\begin{lemma}\cite[Theorem 3.1]{qi2019some}\label{lem:d_kc}
		For $c \in \mathbb{C}, \alpha = \frac{1}{\beta} = \frac{c+\sqrt{c^2-4}}{2}$ and $k \geq 0$, we have
		\begin{align*}
			D_{k}(c) =
			\left\{\begin{aligned}
				\frac{\alpha^{k+1}-\beta^{k+1}}{\alpha-\beta}, & ~~~ c \neq \pm 2;\\
				k + 1, & ~~~ c = 2; \\
				(-1)^k(k+1), & ~~~ c = -2.
			\end{aligned}\right.
		\end{align*}
	\end{lemma}
	Now we shall determine the chain graphs with distinct Seidel eigenvalues.
	\begin{theorem}
		Let $G$ be a chain graph of order $n$ with binary string $0^{a_{1}}1^{a_{2}} \cdots 0^{a_{2h-1}}1^{a_{2h}}$. Then all Seidel eigenvalues of $G$ are distinct if and only if $a_{1} = a_{2} = \cdots = a_{2h-1} = a_{2h} =1.$
	\end{theorem}
	\begin{proof}
		If $G$ has $n$ distinct Seidel eigenvalues, Then we have $h+1 \leq M_{S}(G) = n \leq 2h$ by Reamrk \ref{rem:m_s}. Combining this with $n \geq 2h$, we obtain $n = 2h$, i.e., $a_{1} = a_{2} = \cdots = a_{2h-1} = a_{2h} =1.$\\
		\indent If $a_{1} = a_{2} = \cdots = a_{2h-1} = a_{2h} =1,$ then the Seidel matrix $S(G)$ is the same with the quotient matrix $\widetilde{S}(G)$. Moreover, the multiplicity of the Seidel eigenvalue $-1$ is $1$, and we next consider the multiplicities of other Seidel eigenvalues of $G$. Let $\lambda \neq -1$ be an eigenvalue of $\widetilde{S}(G)$. Replacing $x$ by $\lambda$ in (\ref{equ:matrix_fin}) and replacing $a_{i}$ with $1$ for each $i = 1, \cdots, 2h$, we obtain the matrix $F_{\lambda}$ satisfying $\det(F_{\lambda}) = 0$. In what follows, we shall consider the determinant of the submatrix $F^{(2)}_{\lambda}$ obtained by deleting the last row and last column of $F_{\lambda}$. We will show $\det(F_{\lambda}^{(2)}) \neq 0$ which means that rank($F_{\lambda}) = n-1$ and any Seidel eigenvalue of $G$ is simple. Indeed,
		\begin{gather*}
			\begin{aligned}
				\det(F_{\lambda}^{(2)})& = \begin{vmatrix}
					1+\lambda & 0 & 0 &0 & \cdots & 0 & 0\\
					1+\lambda & 2 & -1-\lambda & 0 & \cdots & 0 & 0\\
					0 & -1-\lambda & 2 & 1+\lambda & \cdots & 0 & 0\\
					0 & 0 & 1+\lambda & 2 & \cdots & 0 & 0\\
					\vdots & \vdots & \vdots & \vdots & \ddots & \vdots & \vdots\\
					0 & 0 & 0 & 0 & \cdots & 2  & -1-\lambda\\
					0 & 0 & 0 & 0 & \cdots & -1-\lambda & 2\\
				\end{vmatrix}_{(2h-1) \times (2h-1)}\\
				& = (1+\lambda)(1+\lambda)^{2h-2}\begin{vmatrix}
					\frac{2}{1+\lambda} & -1 & 0 & 0 & \cdots & 0 & 0\\
					-1 & \frac{2}{1+\lambda} & 1 & 0 & \cdots & 0 & 0\\
					0 & 1 & \frac{2}{1+\lambda} & 1 & \cdots & 0 & 0\\
					\vdots & \vdots & \vdots & \vdots & \ddots & \vdots & \vdots\\
					0 & 0 & 0 & 0 & \cdots & \frac{2}{1+\lambda} & -1\\
					0 & 0 & 0 & 0 & \cdots & -1 & \frac{2}{1+\lambda}
				\end{vmatrix}_{(2h-2) \times (2h-2)}\\
				& = (1+\lambda)^{2h-1}\begin{vmatrix}
					\frac{2}{1+\lambda} & 1 & 0 & 0 & \cdots & 0 & 0\\
					1 & \frac{2}{1+\lambda} & 1 & 0 & \cdots & 0 & 0\\
					0 & 1 & \frac{2}{1+\lambda} & 1 & \cdots & 0 & 0\\
					\vdots & \vdots & \vdots & \vdots & \ddots & \vdots & \vdots\\
					0 & 0 & 0 & 0 & \cdots & \frac{2}{1+\lambda} & 1\\
					0 & 0 & 0 & 0 & \cdots & 1 & \frac{2}{1+\lambda}
				\end{vmatrix}_{(2h-2) \times (2h-2)}\\
				& \neq 0.
			\end{aligned}
		\end{gather*}
		The second equation obtained by the Laplace expansion of $\det(F_{\lambda}^{(2)})$ along with the first row and the basic properties of determinant. The third equation obtained by using the recursive formula of the determinant of any tridiagonal matrix which was introduced in \cite{MEA2004algorithm}. The last inequality follows from Lemma \ref{lem:d_kc}.
	\end{proof}
	
	We also consider the distribution of the Seidel eigenvalues of a chain graph $G$. Keep in mind that $S(G) = J - I - 2A(G)$ and we can use the distribution of $A(G)$ and Courant-Weyl inequalities to derive the following theorem.
	\begin{theorem}\label{tho:se_free}
		Let $G$ be a chain graph with binary string $0^{a_{1}}1^{a_{2}} \cdots 0^{a_{2h-1}}1^{a_{2h}}$, then $G$ has no Seidel eigenvalue in the interval $[-2,0]$. Moreover, the inertia of $S(G)$ is $(n_{+}(S(G)), n_{0}(S(G)),\\ n_{-}(S(G)) = (h, 0, n-h)$.
	\end{theorem}
	\begin{proof}
		Assume that all eigenvalues of $A(G)$ are arranged as follows:
		\begin{gather*}
			\lambda_{1}(A)  \geq \cdots \geq \lambda_{h}(A) \geq \lambda_{h+1}(A) \geq \cdots \geq \lambda_{n-h}(A) \geq \lambda_{n-h+1}(A) \geq \cdots \geq \lambda_{n}(A).
		\end{gather*}
		By Lemma \ref{lem:efree} we know that $A(G)$ has no non-zero eigenvalue in the interval $[-\frac{1}{2}, \frac{1}{2}]$. Moreover, the eigenvalues of $G$ are symmetric with respect to $0$ by the fact that $G$ is bipartite graph. Consequently, we have $\frac{1}{2} < \lambda_{h}(A) = -\lambda_{n-h+1}(A)$  and $\lambda_{h+1}(A) = \cdots = \lambda_{n-h}(A) = 0$.\\
		\indent Similarly, we arrange the eigenvalues of $S(G)$ as follows:
		\begin{gather*}
			\lambda_{1}(S) \geq  \cdots \geq \lambda_{h}(S) \geq \lambda_{h+1}(S) \geq \cdots \geq \lambda_{n-h+1}(S) \geq \lambda_{n-h+2}(S) \geq \cdots \geq \lambda_{n}(S).
		\end{gather*}
		By Lemma \ref{lem:cw} we obtain
		\begin{gather*}
			\lambda_{h}(S) = \lambda_{h}(J-I-2A) \geq \lambda_{n}(J-I) + \lambda_{h}(-2A) = -1 -2\lambda_{n-h+1}(A) > 0,\\
			\lambda_{n-h+2}(S) = \lambda_{n-h+2}(J-I-2A) \leq \lambda_{2}(J-I) + \lambda_{n-h+1}(-2A) = -1 - 2\lambda_{h}(A) < -2.
		\end{gather*}
		Moreover, we have $\lambda_{h+1}(S) = \cdots = \lambda_{n-h+1}(S) = -1$ by the fact that $-1$ is an eigenvalue of $S(G)$ with multiplicity $n-2h+1$.
		This completes the proof.
	\end{proof}
	
	Another problem proposed in \cite{mandal2022spectrum} is the following. By the above theorem, we give a negative answer to this problem.
	\begin{problem}\cite[Problem 2.5]{mandal2022spectrum}
		Is there exist a chain graph $G$ for which the number of positive eigenvalues of $\widetilde{S}(G)$ is not equal to the number of negative eigenvalues?
	\end{problem}
	For any chain graph $G$ with binary string $0^{a_{1}}1^{a_{2}} \cdots 0^{a_{2h-1}}1^{a_{2h}}$, by considering the entries of the spectrum of $S(G)$ and Theorem \ref{tho:se_free}, we conclude that there exists no chain graph such that the number of positive eigenvalues of $\widetilde{S}(G)$ is not equal to the number of negative eigenvalues. However, the problem: ``Is there exist a chain graph $G$ for which the number of distinct positive eigenvalues of $\widetilde{S}(G)$ is not equal to the number of distinct negative eigenvalues?'' have a positive answer. Indeed, both quotient matrices $\widetilde{S}(G_1)$ and $\widetilde{S}(G_2)$ have $2$ distinct positive eigenvalues and $3$ distinct negative eigenvalues, where $G_1$ and $G_2$ are the chain graphs listed in Example \ref{ex:solve_m_s}.\\
	
	\noindent \textbf{Declaration of competing interest}
	
	The authors declare that there is no conflict of interest.
	
	\vskip 0.6 true cm
	\noindent {\textbf{Acknowledgments}}
	
	This research is supported by the National Natural Science Foundation of China (No.11971164).


\baselineskip=0.25in
	
	
	\begin{thebibliography}{99}
		
		\bibitem{andelic2016location} A. Alazemi,  M. An$\dbar$eli\'{c}, S. K. Simi\'{c}, Eigenvalue location for chain graphs, Linear Algebra Appl. 505:194--210, 2016.
		
		\bibitem{alazemi2021eigenvalue} A. Alazemi,  M. An$\dbar$eli\'{c}, T. Koledin, Z. Stani\'{c}, Eigenvalue-free intervals of distance matrices of threshold and chain graphs, Linear Multilinear A. 69 (16): 2959--2975, 2021.
		
		\bibitem{alazemi2022chain}  A. Alazemi, M. An$\dbar$eli\'{c}, K. C. Das, C. M. da Fonseca, Chain graph sequences and Laplacian spectra of chain graphs, Linear Multilinear A. 71 (4): 569--585, 2023.
		
		
		\bibitem{andelic2011on} M. An$\dbar$eli\'{c}, C. M. da Fonseca,  S. K. Simi\'{c}, D. V. To\v{s}i\'{c}, On bounds for the index of double nested graphs, Linear Algebra Appl. 435 (10): 2475--2490, 2011.
		
		\bibitem{andelic2015some} M. An$\dbar$eli\'{c}, E. Andrade, D. M. Cardoso, C. M. da Fonseca, S. K. Simi\'{c}, D. V. To\v{s}i\'{c}, Some new considerations about double nested graphs, Linear Algebra Appl. 483: 323--341, 2015.
		
		
		
		\bibitem{andelic2020tridiagonal} M. An$\dbar$eli\'{c}, C. M. da Fonseca, S.K. Simi\'{c}, Z. Du, Tridiagongal matrices and spectral properites of some graph classes, Czech. Math. J. 70:  1125--1138, 2020.
		
		\bibitem{bell2008GraphsII} F. K. Bell, D. Cvetkovi\'{c}, P. Rowlinson, S. K. Simi\'{c}, Graphs for which least eigenvalue is minimal, II, Linear Algebra Appl. 429(8-9): 2168--2179, 2008.
		
		\bibitem{bhattacharya2008on} A. Bhattacharya, F. Shmuel, U. N. Peled, On the first eigenvalue of bipartite graphs, Electron. J. Combin. 15(R144): 1, 2008.
		
		\bibitem{brouwer2011spectra} A. E. Brouwer, W. H. Haemers, Spectra of Graphs, Springer, New York, 2012.
		
		\bibitem{cvetkovic2010introduction} D. Cvetkovi\'{c}, P. Rowlinson, S. Simi\'{c}, An Introduction to the Theory of Graph Spectra, Cambridge University Press, Cambridge, 2010.
		
		\bibitem{MEA2004algorithm} M. E. EI-Mikkawy, A fast algorithm for evaluating nth order tri-diagonal determinants, J. Appl. Math. Comput. 166: 581--584, 2004.
		
		\bibitem{harary1974which} F. Harary, A. J. Schwenk, Which graphs have integral spectra?, Graphs and Combinatorics, Lecture Notes in Math. 406, Springer-Verlag, Berlin, 45--51, 1974.
		
		\bibitem{mandal2022spectrum} S. Mandal, R. Mehatari, K. C. Das, On the spectrum and energy of Seidel matrix for chain graphs, arXiv preprint arXiv:2205.00310, 2022.
		
		\bibitem{qi2019some} F. Qi, V. Cernanov\'{a}, Y. S. Semenov, Some tridiagonal determinants related to central Delannoy numbers, the Chebyshev polynomials, and the Fibonacci polynomials, Politehn. Univ. Bucharest Sci. Bull. Ser. A Appl. Math. Phys 81 (1): 123--136, 2019
		
		\bibitem{zhuang2022eigenvalue} Z. Xiong, Y. Hou, Eigenvalue-free interval for Seidel matrices of threshold graphs, Appl. Math. Comput. 427: 127177, 2022.
		
	\end{thebibliography}
\end{document}